\documentclass{amsart}

\usepackage{amsmath,amssymb,amsfonts,enumerate,amsthm,amscd,latexsym}

\newcommand{\Max}{\mbox{Max}\,}

\newcommand{\diam}{\mbox{diam}\,}
\newcommand{\clique}{\mbox{clique}\,}

\newcommand{\J}{\mbox{J}}

\newcommand{\G}{\Gamma}
\newcommand{\U}{\mbox{U}}
\newcommand{\I}{\mbox{I}}
\renewcommand{\d}{\mbox{d}}

\newcommand{\fm}{\mathfrak{m}}

\newcommand{\fn}{\mathfrak{n}}
\newtheorem{thm}{Theorem}[section]
\newtheorem{cor}[thm]{Corollary}
\newtheorem{lem}[thm]{Lemma}
\newtheorem{prop}[thm]{Proposition}

\newtheorem{exam}[thm]{Example}

\begin{document}

\bibliographystyle{amsplain}

\author{Hamid Reza Maimani}
\address{Hamid Reza Maimani\\Department of Mathematics, University of
Tehran, Tehran, Iran\\ and Institute for Theoretical Physics and
Mathematics (IPM).}

\email{maimani@ipm.ir}

\author{Maryam Salimi}
\address{Maryam Salimi\\Department of Mathematics, University of
Tehran, Tehran, Iran.}

\author{Asiyeh Sattari}
\address{Asiyeh Sattari\\Department of Mathematics, University of
Tehran, Tehran, Iran.}

\author{Siamak Yassemi}
\address{Siamak Yassemi\\Department of Mathematics, University of
Tehran, Tehran, Iran\\ and Institute for Theoretical Physics and
Mathematics (IPM).}

\email{yassemi@ipm.ir}

\dedicatory{Dedicated to J\"{u}rgen Herzog on the occasion of his
65th birthday}

\thanks{* Corresponding author. Department of Mathematics,
University of Tehran, P.O. Box 13145--448 Tehran, Iran}
\thanks{H. R. Maimani
was supported in part by a grant from IPM No. 85050117}
\thanks{S. Yassemi was supported was supported by a grant from IPM No. 85130214}



\title{Comaximal graph of commutative rings}

\begin{abstract}
Let $R$ be a commutative ring with identity. Let $\Gamma(R)$ be a
graph with vertices as elements of $R$, where two distinct vertices
$a$ and $b$ are adjacent if and only if $Ra+Rb=R$. In this paper we
consider a subgraph $\Gamma_2(R)$ of $\Gamma(R)$ which consists of
non-unit elements. We look at the connectedness and the diameter of
this graph. We completely characterize the diameter of the graph
$\Gamma_2(R)\setminus\J(R)$. In addition, it is shown that for two
finite semi-local rings $R$ and $S$, if $R$ is reduced, then
$\Gamma(R)\cong\Gamma(S)$ if and only if $R\cong S$.

\end{abstract}

\maketitle

\section{Introduction}

For the sake of completeness, first we state some definitions and
notions used throughout to keep this paper as self contained as
possible. We define a {\it coloring} of a graph $G$ to be an
assignment of colors (elements of some set) to the vertices of $G$,
one color to each vertex, so that adjacent vertices are assigned
distinct colors. If $n$ colors are used, then the coloring is
referred to as an {\it $n$-coloring}. If there exists an
$n$-coloring of a graph $G$, then $G$ is called $n$-colorable. The
minimum $n$ for which a graph $G$ is $n$-colorable is called the
{\it chromatic number} of $G$, and is denoted by $\chi(G)$. For a
graph $G$, the {\it degree} of a vertex $v$ in $G$ is the number of
edges of $G$ incident with $v$. Recall that a graph is said to be
{\it connected} if for each pair of distinct vertices $v$ and $w$,
there is a finite sequence of distinct vertices $v=v_1,\cdots,v_n=w$
such that each pair $\{v_i,v_{i+1}\}$ is an edge. Such a sequence is
said to be a path and the distance, $\d(v,w)$, between connected
vertices $v$ and $w$ is the length of the shortest path connecting
them. The {\it diameter} of a connected graph is the supremum of the
distances between vertices. The diameter is 0 if the graph consists
of a single vertex and a connected graph with more than one vertex
has diameter 1 if and only if it is complete; i.e., each pair of
distinct vertices forms an edge. An {\it $r$-partite} graph is one
whose vertex set can be partitioned into $r$ subsets so that no edge
has both ends in any one subset. A {\it complete $r$-partite} graph
is one in which each vertex is joined to every vertex that is not in
the same subset. The {\it complete bipartite} (i.e., $2$-partite)
graph with part sizes $m$ and $n$ is denoted by $K_{m,n}$. A graph
in which each pair of distinct vertices is joined by an edge is
called a {\it complete} graph. We use $K_n$ for the complete graph
with $n$ vertices. A {\it clique} of a graph is its maximal complete
subgraph and the number of vertices in the largest clique of graph
G, denoted by $\clique(G)$, is called the {\it clique number} of
$G$. Obviously $\chi(G)\ge\clique(G)$ for general graph $G$ (see
\cite[page 289]{CO}). Let $G_1=(V_1,E_1)$ and $G_2=(V_2,E_2)$ be two
graphs with disjoint vertices set $V_i$ and edges set $E_i$. The
join of $G_1$ and $G_2$ is denoted by $G=G_1\vee G_2$ with vertices
set $V_1\cup V_2$ and the set of edges is $E_1\cup E_2\cup\{xy|x\in
V_1\,\, \mbox{and}\,\, y\in V_2\}$.

From now on let $R$ be a commutative ring with identity. In
\cite{B}, Beck considered $\G(R)$ as a graph with vertices as
elements of $R$, where two different vertices $a$ and $b$ are
adjacent if and only if $ab=0$. He studied finitely colorable rings
with this graph structure and showed that
$\chi(\G(R))=clique(\G(R))$ for certain classes of rings. in
\cite{AN}, Anderson and Naseer have made further study of finitely
colorable rings and have given an example of a finite local ring
with $5=\clique(\G(R))<\chi(\G(R))=6$.

In \cite{SB}, Sharma and Bhatwadekar define another graph on $R$,
$\Gamma(R)$, with vertices as elements of $R$, where two distinct
vertices $a$ and $b$ are adjacent if and only if $Ra+Rb=R$. They
showed that $\chi(\Gamma(R))<\infty$ if and only if $R$ is a finite
ring. In this case $\chi(\Gamma(R))=\clique(\Gamma(R))=t+\ell$,
where $t$ and $\ell$, respectively, denote the number of maximal
ideals of $R$ and the number of units of $R$.

In this paper, we study further the graph structure defined by
Sharma and Bhatwadekar.

Let $\Gamma_1(R)$ be the subgraph of $\Gamma(R)$, generated by the
units of $R$, and $\Gamma_2(R)$ be the subgraph of $\Gamma(R)$
generated by non-unit elements. In section 2, it is shown that the
graph $\Gamma_2(R)\setminus\J(R)$ is a complete bipartite if and
only if the cardinal number of the set $\Max(R)$ is equal 2 (see
Theorem 2.2). Also we show that $R$ is a finite product of
quasi-local rings if and only if $R$ is clean and
$\clique(\Gamma_2(R)\setminus\J(R))<\infty$ (see Theorem 2.5).

In section 3, the main result says that $\Gamma_2(R)\setminus\J(R)$
is connected and $\diam(\Gamma_2(R)\setminus\J(R))\le 3$ (see
Theorem 3.1). In addition, we completely characterize the diameter
of the graph $\Gamma_2(R)\setminus\J(R)$.

In the final section, it is shown that for two finite semi-local
rings $R$ and $S$, if $R$ is reduced, then $\Gamma(R)\cong\Gamma(S)$
if and only if $R\cong S$ (see Corollary 4.6).

\section{Bipartite graphs}

Throughout this paper $R$ will be a commutative ring with identity,
$\U(R)$ its group of units, $\J(R)$ its Jacobson radical, and
$\I(R)$ its set of idempotents. A ring $R$ is said to be quasi-local
if it has a unique maximal ideal; if $\fm$ is the unique maximal
ideal of $R$, we will often  write $(R,\fm)$.

Let $\Gamma(R)$ be the graph represented by $R$ with definition of
Sharma-Behatwadekar. Let $\Gamma_1(R)=<\U(R)>$ and
$\Gamma_2(R)=<R\setminus\U(R)>$ be the subgraphs of $\Gamma(R)$.
Then it is easy to see that $\Gamma(R)=\Gamma_1(R)\vee\Gamma_2(R)$.

\begin{lem}

The following hold:

\begin{itemize}

\item[(a)] $\Gamma_1(R)$ is a complete graph.

\item[(b)] $a\in\J(R)$ if and only if $\deg_{\Gamma_2(R)}a=0$.

\end{itemize}

\end{lem}

\begin{proof}

Since (a) is clear we just prove (b). Suppose $a\in\J(R)$. Then for
any $\fm\in\Max(R)$, $a\in\fm$. If $\deg_{\Gamma_2(R)}a\neq 0$, then
there exists $b\in\Gamma_2(R)$ such that $Ra+Rb=R$. On the other
hand there exists $\fn\in\Max(R)$ with $b\in\fn$ and so $1\in\fn$
that is a contradiction.

Conversely, assume that $\deg_{\Gamma_2(R)}a=0$. Assume contrary
$a\notin\J(R)$. Then there exists $\fm\in\Max(R)$ such that
$a\notin\fm$. Thus $Ra+\fm=R$. Therefore there exists $b\in\fm$ such
that $Ra+Rb=R$. This contradicts our assumption. In the following we
study the cases where $\Gamma_2(R)\setminus\J(R)$ is complete
bipartite graph and where this graph is $n$-partite.

\end{proof}

We know that each $x\in\U(R)$ is adjacent to every vertex of
$\Gamma(R)$ and it is shown that each $x\in\J(R)$ is an isolated
vertex of $\Gamma_2(R)$. Thus the main part of the graph $\Gamma(R)$
is the subgraph $\Gamma_2(R)\setminus\J(R)$. For this reason the
main aim of this paper is to study the structure of this subgraph.

\begin{thm}

The following are equivalent:

\begin{itemize}

\item[(i)] $\Gamma_2(R)\setminus\J(R)$ is a complete bipartite graph.

\item[(ii)] The cardinal number of the set $\Max(R)$ is equal 2.

\end{itemize}

\end{thm}

\begin{proof}

(ii)$\Rightarrow$(i). Let $\Max(R)=\{\fm_1,\fm_2\}$. Thus the
vertices set of $\Gamma_2(R)\setminus\J(R)$ is equal to the set
$(\fm_1\setminus\fm_2)\cup(\fm_2\setminus\fm_1)$. Let
$a\in\fm_1\setminus\fm_2$ and $b\in\fm_2\setminus\fm_1$. Thus
$Ra+Rb\nsubseteq\fm_1\cup\fm_2$ and so $Ra+Rb=R$.

(i)$\Rightarrow$(ii). Suppose $\Gamma_2(R)\setminus\J(R)$ is a
complete bipartite graph with two part $V_1$ and $V_2$. Set
$M_1=V_1\cup\J(R)$ and $M_2=V_2\cup\J(R)$. We show that $M_1$ and
$M_2$ are two maximal ideals of $R$ and $\Max(R)=\{M_1,M_2\}$. Let
$x,y\in M_1=V_1\cup\J(R)$. Consider the following three cases:

{\it Case 1.} Assume that $x,y\in\J(R)$. Then $x-y\in\J(R)$ and so $x-y\in M_1$.

{\it Case 2.} Assume that $x\in\J(R)$ and $y\in V_1$. Then
$x-y\notin\J(R)$. If $x-y\in\U(R)$, then $Rx+Ry=R$ and so we obtain
a contradiction. If $x-y\in M_2$, then $x-y\in V_2$ and so
$R(x-y)+Ry=R$. Thus $Rx+Ry=R$ which is a contradiction. Therefore
$x-y\in V_1\subseteq M_1$.

{\it Case 3.} Assume that $x,y\in V_1$. If $x-y\in \J(R)$ then there
is nothing to prove. Therefore we assume $x-y\notin\J(R)$. With the
same proof as case 2, the assertion holds.

Now suppose that $r\in R$ and $x\in M_1$. If $x\in\J(R)$, then
clearly $rx\in M_1$. Therefore suppose that $x\notin\J(R)$. Also
$rx$ is not unit. Suppose that $rx\in M_2$. Then $rx\in V_2$ and so
$R(rx)+Rx=R$. Thus $x$ is a unit element of $R$ which is a
contradiction. So $rx\in M_1$.

To now we showed that $M_1$ is an ideal of $R$. By the structure of
$\Gamma(R)$, for any $x\in R\setminus M_1$, we have $M_1+Rx=R$. This
implies that $M_1$ is a maximal ideal.

With the same argument $M_2$ is a maximal ideal of $R$. Now if
$N\in\Max(R)$ then $N\subseteq M_1\cup M_2$ and so $N=M_1$ or
$N=M_2$.\\
This finishes the proof.

\end{proof}

\begin{prop}
Let $n>1$. Then the following hold:

\begin{itemize}

\item[(a)] If $|\Max(R)|=n<\infty$, then the graph $\Gamma_2(R)\setminus\J(R)$ is $n$-partite.

\item[(b)] If the graph $\Gamma_2(R)\setminus\J(R)$ is $n$-partite, then $|\Max(R)|\le n$.
In this case if the graph $\Gamma_2(R)\setminus\J(R)$ is not
$(n-1)$-partite, then $|\Max(R)|=n$.

\end{itemize}

\end{prop}

\begin{proof}

(a). Let $\Max(R)=\{\fm_1,\cdots,\fm_n\}$ and set
$V_1=\fm_1\setminus\J(R)$ and for each $i\ge 2$,
$V_i=\fm_i\setminus\cup_{t=1}^{t=i-1}\fm_t$. Using Prime Avoidence
Theorem, $V_i\neq\varnothing$ for each $i$. It is easy to see that
any two vertices belong to $V_i$ are not adjacent.

(b). Let $V_1,\cdots,V_n$ be the $n$ parts of vertices of
$\Gamma_2(R)\setminus\J(R)$. Assume contrary $|\Max(R)|>n$ and let
$\fm_1,\cdots,\fm_{n+1}\in\Max(R)$. For any $i$, choose
$x_i\in\fm_i\setminus\cup_{j\neq i}\fm_j$. Then it is easy to see
that $\{x_1,\cdots,x_{n+1}\}$ is a clique in
$\Gamma_2(R)\setminus\J(R)$. By the Pigeon Hole Principal, two of
$x_i$'s should belong to one of $V_i$'s, that is a contradiction.
Therefore $|\Max(R)|\le n$. Now suppose that
$\Gamma_2(R)\setminus\J(R)$ is not $(n-1)$-partite and
$|\Max(R)|=m<n$. By (a) the graph will be $m$-partite and this is a
contradiction.

\end{proof}

\begin{prop}

Let $R$ be a ring with $|\Max(R)|\ge 2$. Then the following hold:

\item[(a)] If $\Gamma_2(R)\setminus\J(R)$ is a complete $n$-partite graph, then $n=2$.

\item[(b)] If there exists a vertex of $\Gamma_2(R)\setminus\J(R)$ which is adjacent to every other vertex then $R\cong\mathbb Z_2\times F$, where $F$ is a field.

\end{prop}

\begin{proof} Let $\fm_1, \fm_2$ be two maximal ideals of $R$. Since the elements of $\fm_i\setminus\J(R)$
are not adjacent, and at least one element of $\fm_1\setminus\J(R)$
is adjacent to one element of $\fm_2\setminus\J(R)$, so
$\fm_1\setminus\J(R)$ and $\fm_2\setminus\J(R)$ are subsets of two
distinct parts of $\Gamma_2(R)$. That means
$(\fm_1\setminus\J(R))\cap(\fm_2\setminus\J(R))=\varnothing$. We
claim that $\J(R)=\fm_1\cap\fm_2$. In other case,
$\J(R)\subsetneq(\fm_1\cap\fm_2)$ and so there exists
$x\in(\fm_1\cap\fm_2)\setminus\J(R)$. This elements belongs to
$\fm_1\setminus\J(R)$ and $\fm_2\setminus\J(R)$, that is a
contradiction. Therefore we obtain $\J(R)=\fm_1\cap\fm_2$ and so
$|\Max(R)|=2$. Now by theorem 2.2 we have $n=2$.

(b). Let $x$ be a non-unit element of $R$ which is adjacent to every
other vertex of $\Gamma_2(R)\setminus\J(R)$. Since $x$ is comaximal
with each nonunit outside the Jacobson radical, $x$ is idempotent,
$\J(R)=(0)$ and $\fm=\{0,x\}$ is a maximal ideal. Thus for each
nonunit $s\in R\setminus\fm$, having $xR+sR=R$ implies
$(1-x)sR=(1-x)R$ and this implies $(1-x)R=F$ is a field. Hence
$R\cong{\Bbb Z}_2\times F$.


\end{proof}

A ring is said to be {\it clean} if each of its elements can be
written as the sum of a unit and an idempotent cf. \cite{N} (see
also \cite{AC}). For example, a quasi-local ring is clean. The
following result gives an application of Sharma-Bhatwadegar graph to
characterize clean rings.

\begin{thm}

For the ring $R$, the following are equivalent:

\begin{itemize}

\item[(a)] $R$ is a finite product of quasi-local rings.

\item[(b)] $R$ is clean and $\clique(\Gamma_2(R)\setminus\J(R))$ is finite.

\end{itemize}

\end{thm}

\begin{proof}

(a)$\Rightarrow$(b). Let $R=R_1\times\cdots\times R_n$ where each
$R_i$ is quasi-local with unique maximal ideal $\fm_i$. Set
$N_i=R_1\times\cdots\times R_{i-1}\times\fm_i\times
R_{i+1}\times\cdots\times R_n$ for any $i=1,\cdots,n$. Then each
$N_i$ belongs to $\Max(R)$. For any $i$ choose $x_i\in
N_i\setminus\cup_{\ell\neq i}N_{\ell}$. Then it is easy to see that
$Rx_i+Rx_j=R$ for all $i\neq j$. In addition by using the Pigeon
Hole Principal, there is no any $n+1$ family elements of
$\cup_{i=1}^nN_i$ which pairwise adjacent. Thus
$\clique(\Gamma_2(R)\setminus\J(R))=n<\infty$.\\
On the other hand, each $R_i$ is clean and so by \cite[Proposition
2(3)]{AC}, $R$ is clean.

(b)$\Rightarrow$(a). Suppose that
$\clique(\Gamma_2(R)\setminus\J(R))$ is finite. Assume contrary that
$\I(R)$ has infinitely many idempotent elements then by \cite[Lemma
2.1]{SB} there exists an infinite sequence $e_1,e_2,\cdots$ of
non-trivial idempotents in $\Gamma_2(R)\setminus\J(R)$ such that the
set $S$ consisting of elements $e_i$ ($i\ge 1$) is an infinite
clique. This is a contradiction.

\end{proof}

\section{Diameter of the graph}

In this section we completely characterize the diameter of
$\Gamma_2(R)\setminus\J(R)$. The following result shows that
$\Gamma_2(R)\setminus\J(R)$ is a connected graph and its diameter is
not greater than 3.

\begin{thm}
The graph $\Gamma_2(R)\setminus\J(R)$ is connected, and $\diam(\Gamma_2(R)\setminus\J(R))\le 3$.

\end{thm}

\begin{proof}
Let $a,b\in(R\setminus\U(R))\setminus\J(R)$. We consider two cases:

{\it Case 1} Assume that $ab\notin\J(R)$. There exists
$x\in(R\setminus\U(R))\setminus\J(R)$ such that $Rab+Rx=R$. Thus
$Ra+Rx=Rb+Rx=R$. So we have the path $a$---$x$---$b$, and so
$\d(a,b)\le 2$.

{\it Case 2} Assume that $ab\in\J(R)$. Set $S_a=\{\fm|\fm\in\Max(R),
a\in\fm\}$ and $S_b=\{\fm|\fm\in\Max(R), b\in\fm\}$. Clearly,
$\Max(R)=S_a\cup S_b$. Now suppose that $x$ is adjacent to $a$ in
$\Gamma_2(R)$. Then $x\notin\J(R)$. If $a\in\fm$, then $x\notin\fm$
and so $x\in\fn\in\Max(R)$, where $\fn\in S_b\setminus S_a$. Thus
$bx\notin\J(R)$. Therefore by Case 1, $\d(b,x)\le 2$ and so
$\d(a,b)\le 3$.

\end{proof}

\begin{lem}
$\diam(\Gamma_2(R)\setminus\J(R))=1$ if and only if $R\cong\mathbb
Z_2\times\mathbb Z_2$.

\end{lem}

\begin{proof}
If $\diam(\Gamma_2(R)\setminus\J(R))=1$, then
$\Gamma_2(R)\setminus\J(R)$ is complete graph. Thus there exists a
vertex of $\Gamma_2(R)\setminus\J(R)$ which is adjacent to every
other vertex. Therefore $R\cong\mathbb Z_2\times F$, where $F$ is a
field by Proposition 2.4(b). Since $\Gamma_2(R)\setminus\J(R)$ is
complete, we have that $F\cong\mathbb Z_2$. Thus $R\cong\mathbb
Z_2\times\mathbb Z_2$.

It is easy to see that for $R\cong\mathbb Z_2\times\mathbb Z_2$,
$\diam(\Gamma_2(R)\setminus\J(R))=1$.

\end{proof}

Our next result characterizes the graphs where $\diam(\Gamma_2(R)\setminus\J(R))=2$.

\begin{prop}
Assume that $R$ is not local. The diameter of the graph
$\Gamma_2(R)\setminus\J(R)$ is equal 2 if and only if one of the
following holds:

\begin{itemize}

\item[(a)] $\J(R)$ is a prime ideal.

\item[(b)] $|\Max(R)|=2$ and $R\ncong\mathbb Z_2\times\mathbb Z_2$.

\end{itemize}

\end{prop}

\begin{proof}

Note that if $\J(R)$ is prime and $R$ is semi-local (i.e. has finite
number of maximal ideals), then $R$ will be local. Let $\J(R)$ be a
prime ideal and $a,b\notin\J(R)$. Then $ab\notin\J(R)$, and so by
the same argument as Theorem 3.1, there exists
$x\in(R\setminus\U(R))\setminus\J(R)$, such that $a$---$x$---$b$ is
a path. Thus $\diam(\Gamma_2(R)\setminus\J(R))\le 2$. If
$\diam(\Gamma_2(R)\setminus\J(R))=1$, then by previous result
$R\cong\mathbb Z_2\times\mathbb Z_2$. But $\J(\mathbb
Z_2\times\mathbb Z_2)$ is not a prime ideal. That is a
contradiction.

Now let $|\Max(R)|=2$ and $R\ncong\mathbb Z_2\times\mathbb Z_2$,
then by Theorem 2.2, $\Gamma_2(R)\setminus\J(R)$ is a complete
bipartite graph where at least one of the parts has at least two
elements. Therefore $\diam(\Gamma_2(R)\setminus\J(R))=2$.

Conversely, suppose that $\diam(\Gamma_2(R)\setminus\J(R))=2$ and
$\J(R)$ is not prime. let $a,b\notin\J(R)$ but $ab\in\J(R)$. We
claim that $a$ and $b$ are adjacent. Otherwise, there exists $t$ in
$\Gamma_2(R)$ such that $Ra+Rt=Rb+Rt=R$. Thus $Rab+Rt=R$ and so
$ab\notin\J(R)$ which is a contradiction. Therefore $Ra+Rb=R$ and so
for some $p,q\in R$, $pa+qb=1$. Set $S=R/\J(R)$ and $a_1=pa+\J(R)$
and $b_1=qb+\J(R)$. Then $a_1b_1=0$ and $a_1+b_1=1_S$. Therefore
$a_1$ and $b_1$ are idempotent elements in $S$, and so $S=Sa_1\oplus
Sb_1$. We will show that $Sa_1$ is a field. Let $0\neq x\in Sa_1$
and $0\neq y\in Sb_1$. Then there exists $\alpha, \beta$ such that
$\alpha x+\beta y=1_S$ and so
$\alpha(a_1+b_1)x+\beta(a_1+b_1)y=1_S$. Thus $(\alpha a_1)x+(\beta
b_1)y=1_S$. On the other hand $a_1+b_1=1_S$ and so $(\alpha
a_1)x=a_1$ and $(\beta b_1)y=b_1$. Therefore $x$ is a unit in
$Sa_1$. Therefore $Sa_1$ and $Sb_1$ are fields and so $|\Max(S)|=2$.
therefore $|\Max(R)|=2$.

\end{proof}

\begin{exam}
Let $R=\mathbb Z_n$ where $n=p_1^{\ell_1}\cdots p_r^{\ell_r}$.

Assume $r\ge 3$. Let $x=p_1^{\ell_1}\cdots p_{r-1}^{\ell_{r-1}}$ and
$y=p_2^{\ell_2}\cdots p_{r}^{\ell_{r}}$. Then $x$ and $y$ are not
adjacent. Also if $x,y$ are adjacent $z$, then $(z,x)=(z,y)=1$,
which is impossible. We have
$Rx+Rp_r^{\ell_r}=R=Rp_r^{\ell_r}+Rp_1^{\ell_1}=Rp_1^{\ell_1}+Ry$.
Hence there is path $x$---$p_r^{\ell_r}$---$p_1^{\ell_1}$---$y$. So
$\diam(\Gamma_2(\mathbb Z_n)\setminus\J(\mathbb Z_n))=3$.

Assume that $r=2$. In this case we have two maximal ideals
$M_1=<p_1>$ and $M_2=<p_2>$. Then $\Gamma(R)$ is a complete
bipartite graph and so $\diam(\Gamma_2(\mathbb
Z_n)\setminus\J(\mathbb Z_n))=2$.

Assume that $r=1$. Then $R$ is local and so $\Gamma_2(\mathbb Z_n)\setminus\J(\mathbb Z_n)$ is empty graph.

\end{exam}

\begin{exam}

Let $R$ be an infinite PID. Then for any two non-unit elements
$a,b$, there exists a prime element $p$ such that $p$ does not
divide $a$ and $b$. Therefore $Ra+Rp=Rb+Rp=1$. So $\d(a,b)\le 2$ and
hence $\diam(\Gamma_2(R)\setminus\J(R))=2$.

\end{exam}

\section{isomorphisms}

Recall that two graphs $G$ and $H$ are isomorphic, denoted by
$G\cong H$, if there is a bijection $\varphi : G\to H$  of vertices
such that the vertices $x$ and $y$ are adjacent in $G$ if and only
if $\varphi(x)$ and $\varphi(y)$ are adjacent in $H$.

In this section, we consider the following question:

If $R$ and $S$ are two rings with $\Gamma(R)\cong\Gamma(S)$, then do
we have $R\cong S$?

The following examples show that the above question is not valid in
general.

\begin{exam} Let $R=\mathbb Z_4$ and $S=\mathbb Z_2[x]/(x^2)$. Then by simple computation we can see
that $\Gamma(R)\cong\Gamma(S)(\cong K_2\vee\bar{K_2})$. But $\mathbb Z_4$ and $\mathbb Z_2[x]/(x^2)$ are not isomorphic.
\end{exam}

\begin{exam}
Let $R=\mathbb Z_8$ and $S=\mathbb Z_2[x]/(x^3)$. Then
$\Gamma(R)\cong\Gamma(S) (\cong K_4\vee\bar{K_4})$. But $R$ and $S$
are not isomorphic.

\end{exam}

\begin{exam}

Let $R=\mathbb Z_2[x]/(x^3)$ and $S=\mathbb Z_2[x,y]/(x^2,y^2,xy)$.
Then $\Gamma(R)\cong\Gamma(S) (\cong K_4\vee \bar{K_4})$). But $R$
and $S$ are not isomorphic.

\end{exam}

In the following theorem we give a partial answer to the above
question.

\begin{thm} Let $\{(R_i,\fm_i)\}_{i=1}^m$ and $\{(S_j,\fn_j)\}_{j=1}^n$ be two finite
families of finite quasi-local rings, and let
$R=R_1\times\cdots\times R_m$ and $S=S_1\times\cdots\times S_n$. If
$\Gamma(R)\cong\Gamma(S)$ then $m=n$ and there is a permutation
$\sigma$ on the set $\{1,2,\ldots ,m\}$ such that
$|R_i/\fm_i|=|S_{\sigma(i)}/\fn_{\sigma(i)}|$ for each $i=1,\cdots,
m$, and hence $R_i/\fm_i\cong S_{\sigma(i)}/\fn_{\sigma(i)}$. In
particular, if $\Gamma(R)\cong\Gamma(S)$ and each $R_i$ is a finite
field, then each $S_j$ is also a finite field and $R_i\cong
S_{\sigma(i)}$ for each $i=1,\cdots, m$, and thus $R\cong S$.

\end{thm}

\begin{proof} First note that since $|\Max(R)|=n$ and $|\Max(S)|=m$, and $\Gamma(R)\cong\Gamma(S)$,
we have that $m=n$. Set $M_i=R_1\times\cdots\times
R_{i-1}\times\fm_i\times R_{i+1}\times\cdots\times R_m$ and
$N_i=S_1\times\cdots\times S_{i-1}\times\fn_i\times
S_{i+1}\times\cdots\times S_m$ for each $i=1,\cdots,m$. For any
$i=1,\cdots ,m$ let $x_i\in M_i\setminus\cup_{j\neq i}M_j$. Clearly
$\{x_1,\cdots, x_m\}$ is a clique in $\Gamma_2(R)$. Suppose that
$\nu(x_i)$ is equal to the number of vertices of $\Gamma(R)$ which
are not adjacent to $x_i$. Then $$\nu(x_i)=|R-\{y|y=x_i\,\,\mbox{or
$y$ adjacent to $x_i$}\}|=|M_i|.$$ For each $i=1,\cdots ,m$, let
$y_{\sigma(i)}\in N_{\sigma(i)}$ be the image of $x_i$ under the
graph isomorphism. Then $\{y_1,\cdots,y_m\}$ is a clique in
$\Gamma_2(S)$ and $y_{\sigma(i)}\in
N_{\sigma(i)}\setminus\cup_{j\neq \sigma(i)}N_j$. It is easy to see
that $\nu(y_{\sigma(i)})=|N_{\sigma(i)}|$. Thus
$|M_i|=|N_{\sigma(i)}|$ and so
$|R_i/\fm_i|=|S_{\sigma(i)}/\fn_{\sigma(i)}|$. Therefore
$R/\fm_i\cong S_{\sigma(i)}/\fn_{\sigma(i)}$.

In particular, if $\Gamma(R)\cong\Gamma(S)$ and each $R_i$ is a
finite field. Thus $\J(R)=(0)$ and so $\J(S)=(0)$ and hence
$\fn_i=(0)$ for each $i$. Therefore each $S_j$ is also a finite
field and $R_i\cong S_{\sigma(i)}$ for each $i\in I$, and thus
$R\cong S$.

\end{proof}

The following example shows that the condition ``$R_i$ is a field''
is necessary in Theorem 4.4.

\begin{exam}
Let $R=\mathbb Z_2\times \mathbb Z_8$, and $S=\mathbb
Z_4\times\mathbb Z_4$. Then $\J(R)\cong\{0\}\times\mathbb Z_4$ and
$\J(S)\cong\mathbb Z_2\times\mathbb Z_2$. Also
$\U(R)=\{1\}\times\{1,3,5,7\}$ and $\U(S)=\{1,3\}\times\{1,3\}$.
Since $|\Max(R)|=|\Max(S)|=2$, then
$\Gamma_2(R)\setminus\J(R)\cong\Gamma_2(S)\setminus\J(S)\cong
K_{4,4}$. Therefore $\Gamma(R)\cong\Gamma(S)\cong(K_{4,4}\cup
\bar{K_4})\vee K_4$. But it is clear that $R\ncong S$.

\end{exam}

\begin{cor} Let $R$ and $S$ be two finite semi-local rings and let $R$ be reduced.
Then $\Gamma(R)\cong\Gamma(S)$ if and only if $R\cong S$.

\end{cor}

\begin{proof} It is clear that $R=F_1\times\cdots\times F_n$, where $F_i$ is a field for any $i=1,\cdots,n$.
Now the assertion holds from Theorem 4.4.

\end{proof}

The following result shows that there exists a copy of
$\Gamma(R/\J(R))$ in the structure of $\Gamma(R)$. This result
obtains that for two rings $R$ and $S$ if $\Gamma(R)\cong\Gamma(S)$,
then $R/\J(R)\cong S/\J(S)$.

\begin{prop}
The following hold:

\begin{itemize}

\item[(a)] If $a$ is adjacent to $b$ in $\Gamma(R)$, then every element of $a+\J(R)$ is adjacent to $b+\J(R)$.

\item[(b)] The elements of $a+\J(R)$ are adjacent if and only if $a$ is an unit. In this case,
each element of $a+\J(R)$ is unit too.

\item[(c)] There exists a copy of $\Gamma(R/\J(R))$ in the structure of $\Gamma(R)$.
In particular, if $\Gamma(R)\cong\Gamma(S)$, then $R/\J(R)\cong S/\J(S).$

\end{itemize}

\end{prop}

\begin{proof}

(a). Suppose that $Ra+Rb=R$. Let $x=a+r_1$ and $y=b+r_2$ for $r_1,r_2\in\J(R)$.
Then there exists elements $s,t\in R$ such that $sa+tb=1$. So
\[ \begin{array}{rl}
 sx+ty &\, =sa+tb+sr_1+tr_2\\
 &\, = 1-(-sr_1-tr_2).
\end{array} \]
Since $-sr_1-tr_2\in\J(R)$, so $sx+ty$ is a unit and hence $Rx+Ry=R$.

(b). Let $(a+r_1)$ be adjacent to $a+r_2$. Then
$R(a+r_1)+R(a+r_2)=R$ and so there exist $s,t\in R$ such that
$s(a+r_1)+t(a+r_2)=1$. This implies that $(s+t)a=1-(r_1s+r_2t)$ and
so $(s+t)a$ is invertible. Therefore $a$ is invertible.

(c). Choose a distinct representation $\{a_i\}$ from the cosets of
$R/\J(R)$. By parts (a) and (b), we have
$<\{a_i\}>\cong\Gamma(R/\J(R))$. Let $\varphi:\Gamma(R)\to\Gamma(S)$
be an isomorphism. Then $\varphi(<\{a_i\}>)=\{b_i\}$ and
$<\{b_i\}>=S/\J(S)$ and hence the assertion is easily obtained.
\end{proof}

\section*{Acknowledgments}

This paper was finalized when H.R. Maimani and S. Yassemi were
visiting the Tata Institute of Fundamental Research (TIFR) under
TWAS-UNESCO Associateship Scheme. It is a pleasure to thank both
TWAS and TIFR for financial support and hospitality. The authors
would like to thank S.M. Bhatwadekar for the stimulating
discussions. The authors wish to thank an anonymous referee, whose
comments have improved this paper.

\providecommand{\bysame}{\leavevmode\hbox
to3em{\hrulefill}\thinspace}

\end{document}